\title{Regular induced subgraphs of a random graph}
\author{
Michael Krivelevich\thanks{School of Mathematical Sciences, Tel Aviv University, Tel Aviv 69978, Israel.
E-mail: {\tt krivelev@post.tau.ac.il}. Research supported
  in part by a USA-Israel BSF Grant, by a grant from the
  Israel Science Foundation, and by Pazy Memorial Award.}
  \and Benny Sudakov\thanks{Department of Mathematics, UCLA 90095, Los Angeles, USA.
  E-mail: {\tt bsudakov@math.ucla.edu}. Research supported in part by NSF CAREER award DMS-0546523 and a USA-Israeli BSF grant.}
 \and  Nicholas  Wormald \thanks{Department of Combinatorics and Optimization, University of Waterloo,  Waterloo ON, Canada. E-mail: 
{\tt  nwormald@uwaterloo.ca}.    Supported by the  Canada Research Chairs Program and NSERC.}}

\documentclass[11pt]{article}

\usepackage{amsmath}
\usepackage{amsfonts}
\usepackage{latexsym}

  \newcommand{\lab}[1]{\label{#1}}                

\newtheorem{theorem}{Theorem}[section]
\newtheorem{lemma}[theorem]{Lemma}
\newtheorem{corollary}[theorem]{Corollary}

\newtheorem{prop}[theorem]{Proposition}

\newcommand{\kh}{(k-1)/2}
\newcommand{\whp}{\emph{whp}}
\date{}

\def\proof{\noindent{\bf Proof.}\ \ }
\def\qed{\hfill\vrule height8pt width4pt depth0pt}
\newcommand{\bel}[1]{\be\lab{#1}}
\newcommand\eqn[1]{(\ref{#1})}
\newcommand{\be}{\begin{equation}}
\newcommand{\ee}{\end{equation}}
\newcommand{\bea}{\begin{eqnarray}}
\newcommand{\eea}{\end{eqnarray}}
\newcommand{\bean}{\begin{eqnarray*}}
\newcommand{\eean}{\end{eqnarray*}}

\def\kk{{k}}
\def\pr{{\mathbb{P}}}
\def\ex{{\mathbb{E}}}
\def\eps{{\epsilon}}

\def\la{{\lambda}}
\def\dv{{\bf d}}
\def\cv{{\bf c}}
\def\Dscr{{\cal D}}
\def\dbar{{\bar d}}
\def\dhat{{\hat d}}
\def\sv{{\bf s}}
\def\Bscr{{\cal B}}
\def\Sscr{{\cal S}}
\def\Wscr{{\cal W}}

\begin{document}
\maketitle
\begin{abstract}
An old problem of Erd\H{o}s, Fajtlowicz and Staton asks for the order of a largest induced regular subgraph that can be found in every graph on $n$ vertices. Motivated by this problem, we consider the order of such a subgraph in a typical graph on $n$ vertices, i.e., in a binomial random graph $G(n,1/2)$. We prove that with high probability a largest induced regular subgraph of $G(n,1/2)$ has about $n^{2/3}$ vertices.
\end{abstract}

\section{Introduction}
A rather old and apparently quite difficult problem of  Erd\H{o}s,
Fajtlowicz and Staton (see \cite{Er1} or \cite{CG}, page 85)
asks for the order of a largest induced
regular subgraph that can be found in every graph on $n$ vertices.
By the known estimates for graph Ramsey numbers (c.f., e.g.,
\cite{GRS}), every graph on $n$ vertices contains a clique
or an independent set of size $c\ln n$, for some positive constant
$c>0$, providing a trivial lower bound of $c\ln n$ for the problem.
Erd\H{o}s, Fajtlowicz and Staton conjectured that the quantity in
question is asymptotically larger than $\log n$. 
So far this conjecture has not been settled. Some progress has been
achieved in upper bounding this function of $n$: Bollob\'as in an
unpublished argument 
showed (as stated in~\cite{CG}) the existence of a graph on $n$ vertices without an induced
regular subgraph on at least $n^{1/2+\epsilon}$ vertices, for any
fixed $\epsilon>0$ and sufficiently large $n$. A slight
improvement has recently been obtained by Alon and the first two
authors \cite{AKS}, who took the upper bound down to
$cn^{1/2}\log^{3/4}n$. 

Given the simplicity of the problem's statement, its appealing
character and apparent notorious difficulty, it is quite natural
to try and analyze the behavior of this graph theoretic parameter
for a typical graph on $n$ vertices, i.e. a graph drawn from the
probability space $G(n,1/2)$ of graphs. (Recall that the ground
set of the probability space $G(n,p)$ is composed of all graphs on
$n$ labeled vertices, where each pair $(i,j)$ appears as an edge
in $G$, drawn from $G(n,p)$, independently and with probability $p$.
In the case $p=1/2$ all labeled graphs $G$ on $n$ vertices are
equiprobable: $Pr[G]=2^{-{n\choose 2}}$.) This is the subject of
the present paper. 

 We say that a graph property ${\cal P}$ holds with high probability, or \whp\ for brevity, if the probability of a random graph to have ${\cal P}$ tends to 1 as $n$ tends to infinity.  We prove the following result.
\begin{theorem}\label{th1}
Let $G$ be a random graph $G(n,1/2)$. Then with high probability every induced regular subgraph of $G$ has at most $2n^{2/3}$ vertices. On the other hand, for  $k=o(n^{2/3})$, with high probability
$G$ contains a set of $k$ vertices that span  a $(k-1)/2$-regular graph.
\end{theorem}

It is instructive to compare this result with the above
mentioned result of \cite {AKS}. Alon et al.\ also used a certain
probability space of graphs to derive their upper bound of
$O(n^{1/2}\log^{3/4}n)$. Yet, their model of random graphs is much
more heterogeneous in nature (the expected degrees of vertices
vary significantly there, see \cite{AKS} for full details). As
expected, the rather homogeneous model $G(n,1/2)$ produces a
sizably weaker upper bound for the  Erd\H{o}s-Fajtlowicz-Staton
problem. 

The difficult part of our proof is the lower bound. For this we use the second moment method. Getting an accurate bound on the variance is the main difficulty. Our main tool for this bounds the number of regular 
graphs on $k$ vertices which contain given subgraph $H$, when $H$ is not too large. For $H$ with $o(\sqrt k)$ vertices and with degree sequence satisfying certain conditions, we obtain an asymptotic formula for this number which is of independent interest; see Theorem~\ref{c:probformula}. 

In Section 2 we introduce some notation and technical tools
utilized in our arguments, and then prove a rather straightforward
upper bound in Theorem \ref{th1}. A much more delicate lower bound
is then proven in Section~3. The technical lemma used in this proof relies on the above-mentioned estimate of the number of regular graphs with a given subgraph. Its proof is relegated to Section~4. Section~5, the final section of the paper, contains some
concluding remarks.

\section{Notation, tools and the upper bound}
In this short section we describe some notation and basic tools to be used later in our proofs. Then we establish the upper bound part of Theorem \ref{th1}.

We will utilize the following (standard) asymptotic notation.
For two functions $f(n)$, $g(n)$ of
a natural valued parameter $n$, we write $f(n)=o(g(n))$, whenever
$\lim_{n \rightarrow\infty} f(n)/g(n)=0$; $f(n)=\omega(g(n))$ if $g(n)=o(f(n))$.
Also, $f(n)=O(g(n))$ if there
exists a constant $C>0$ such that $f(n)\le C g(n)$ for all $n$;
$f(n)=\Omega(g(n))$ if $g(n)=O(f(n))$, and $f(n)=\Theta(g(n))$
if both $f(n)=O(g(n))$ and $f(n)=\Omega(g(n))$ are satisfied. 
We write $f\sim g$ if the ratio $f/g$ tends to 1 when the underlying parameter tends to infinity.
For a real $x$ and positive integer $a$, define $[x]_a=x(x-1)\cdots(x-a+1)$.  
All logarithms in this paper have the natural basis. We will use the bound ${n\choose k}=(en/k)^k$, valid for all positive $n$ and $k$.

Let $G(\dv)$ denote the number of labeled simple graphs on
$k$ vertices with degree
sequence $\dv =(d_1,d_2,\ldots,d_k)$. Also, we denote
$$
p_k = \mathbb{P} [\mbox{a random graph $G(k,0.5)$ is $\lfloor\kh\rfloor$-regular}]\,.
$$
Clearly, $p_k=G(\dv)2^{-{k\choose 2}}$, with all $d_i$ being equal to $\lfloor\kh\rfloor$.

We will cite repeatedly the following corollary of a result of McKay and the third author (see Theorems 2 and 3 of
\cite{MWhigh}).
\begin{theorem}\lab{t:counting}
Let $d_j=d_j(k)$, $1\le j\le k$ be integers such that
$\sum_{j=1}^k d_j=\la k(k-1)$ is an even integer where $1/3<\la<2/3$,
and $|\la k-d_j|=O(k^{1/2+\eps})$ uniformly over $j$, for some sufficiently small fixed $\eps$.  Then
\bel{formula}
G(\dv) =f(\dv) \big(\la^\la (1-\la)^{1-\la}\big)^{k\choose 2}
        \prod_{j=1}^k {{k-1}\choose d_j}
\ee
where
\begin{itemize}
\item $f(\dv)=O(1)$, and
\item if $\max\{|\la k-d_j|\}=o(\sqrt{k})$, then $f(\dv)\sim \sqrt{2}e^{1/4}$, uniformly over the choice of such a
degree sequence $\dv$.
\end{itemize}
 \end{theorem}
Observe that the expression  $\big(\la^\la (1-\la)^{1-\la}\big)^{k\choose 2}
\prod_{j=1}^k {{k-1}\choose d_j}$ is at most $2^{-{k\choose 2}}{{k-1}\choose {\lfloor\kh\rfloor}}^k$, hence $G(\dv)2^{-{k\choose 2}}$ is $O(p_k)$ for every degree sequence $\dv$ covered by Theorem \ref{t:counting}. Also,
using Stirling's formula is it straightforward to verify that $p_k=\big((1+o(1))\sqrt{\pi k/2}\big)^{-k}$ and that $p_{k-1}/p_k=\Theta(\sqrt{k})$.

In order to prove the upper bound in Theorem~\ref{th1}, we show that, for a given $k$ and $r$, the probability that a random graph on $k$ vertices is $r$-regular is
$O(p_k)$.  (For future use we prove here a somewhat more general statement.) We then use the  above-mentioned  estimate for  $p_k$ and apply the union bound over all possible values of $r$.

\begin{lemma}\label{newclaim1}
For every degree sequence $\dv=(d_1,\ldots,d_k)$,
$$
\mathbb{P}[G(k,0.5)\mbox{ has degree sequence $\dv$}]= O(p_k)\ .
$$
\end{lemma}

\proof
Let $\dv$ be a degree sequence of length $k$ for which $G(\dv)$ is maximal (which is obviously equivalent to choosing $\dv$ to be a most probable degree sequence in $G(k,1/2)$).
If all degrees in $\dv$ satisfy
$|d_i-k/2|\le k^{1/2+\epsilon}$, then Theorem \ref{t:counting} is applicable, and we are done. Otherwise, there is
$d_i$, say, $d_k$, deviating from $k/2$ by at least
$k^{1/2+\epsilon}$, for some fixed $\epsilon>0$. To bound the probability that $G(k,1/2)$ has degree sequence $\dv$, we first expose the edges from vertex $k$ to the rest of the graph.  By standard estimates on the tails of the binomial distribution, the probability that $k$ has the required
degree is $\exp\{-\Omega(k^{2\epsilon})\}$. The edges exposed induce a
new degree sequence on vertices $1,\ldots,k-1$. Observe that in order to contradict the lemma's assertion there should be some degree sequence $\dv'$ of length $k-1$, whose probability in $G(k-1,1/2)$ is larger than $p_k$ by the exponential factor of  $\exp\{\Omega(k^{2\epsilon})\}$. Since the ratio $p_{k-1}/p_k$ is of order $\Theta(\sqrt{k})$, it follows that the probability of $\dv'$ to appear in $G(k-1,1/2)$ is at least
$p_{k-1}\cdot\exp\{\Omega(k^{2\epsilon})\}$. Repeating this argument at most $k/2$ times we either prove the lemma or conclude that there should exist a degree sequence $\dv''$ of length $k/2$ whose probability in $G(k/2,1/2)$ is at least $p_{k/2}\cdot\exp\{\Omega(k^{2\epsilon})\cdot k/2\}$. Recalling that $p_{k/2}=\big(\Theta(1/\sqrt{k})\big)^{k/2}$, the latter expression is more than 1 ---     a contradiction.
\qed

\medskip

In order to complete the proof of the upper bound of Theorem \ref{th1}, note that by Lemma \ref{newclaim1} the probability that a fixed set $V_0$ of $k$ vertices spans a regular subgraph in $G(n,1/2)$ is $O(kp_k)$. Summing over all $k\ge k_0=2n^{2/3}$ and all vertex subsets of size $k$, we conclude that the probability that $G(n,1/2)$ contains an induced regular subgraph on at least $k_0$ vertices is 
$$
\sum_{k\ge k_0} {n\choose k}\cdot O(kp_k) \le\sum_{k\ge k_0} \left(\frac{en}{k}\right)^k\,k \big((1+o(1))\sqrt{\pi k/2}\big)^{-k}
\le n^2\cdot \left(\frac{(1+o(1))\sqrt{2}en}{\sqrt{\pi}k_0^{3/2}}\right)^{k_0}= o(1)\,.
$$

\section{A lower bound}

In this section we give a proof of the lower bound in our main result, Theorem \ref{th1}. (To be more accurate, we give here most of the proof, deferring the proof of a key technical lemma to the next section.) The proof uses the so-called second moment method and proceeds by estimating carefully the first two moments of the random variable $X=X(k)$, counting the number of $\kh$-regular induced subgraphs on $k$ vertices in $G(n,1/2)$.
For convenience we assume throughout the proof that $k$ is odd. We find it quite surprising that it is possible to apply the second moment method to sets of such a large size.

So let $X$ be the random variable counting the number of
$\kh$-regular induced subgraphs on $k$ vertices in $G(n,0.5)$.
We write $X=\sum_{|A|=k} X_A$, where $X_A$ is the indicator random variable for the event that a vertex subset $A$ spans a $\kh$-regular subgraph.
Then
$$
\mathbb{E}[X] =\sum_{|A|=k} \mathbb{E}[X_A]={n\choose k} p_k\,.
$$
Plugging in the estimate for $p_k$ cited after the statement of Theorem  \ref{t:counting}, it is straightforward to verify that $\mathbb{E}[X]$ tends to infinity for $k=o(n^{2/3})$; in fact, $\mathbb{E}[X]=(\omega(1))^k$ in this regime. Denote by $\mathbb{V}ar[X]$ the variance of $X$. A corollary of Chebyshev's inequality is that $\mathbb{P}[X>0]\ge 1-\frac{\mathbb{V}ar[X]}{\mathbb{E}^2[X]}$, and therefore in order to prove that \whp\ $G(n,1/2)$ contains an induced regular subgraph on $k$ vertices, it is enough to establish that $\mathbb{V}ar[X]=o(\mathbb{E}^2[X])$.

In order to estimate the variance of $X$ we need to estimate the correlation between the following events: ``$A$ spans a $\kh$-regular subgraph" and  ``$B$ spans a $\kh$-regular subgraph", where $A,B$ are $k$-element vertex subsets whose intersection is of size $i\ge 2$. To this end, define
$$
p_{k,i} = \max_{|H|=i}\mathbb{P}[\mbox{$G(k,0.5)$ is $\kh$-regular} \mid
G[i]=H] ,
$$
where the maximum in the expression above is taken over all graphs $H$ on $i$ vertices, and $G[i]$ stands for the subgraph of $G(k,1/2)$ spanned by the first $i$ vertices.
Since $X=\sum_{|A|=k} X_A$, we have:
\begin{eqnarray}
\mathbb{V}ar[X]&=&\mathbb{E}[X^2]-\mathbb{E}^2[X]= \sum_{|A|=k} \mathbb{V}ar[X_A]+
\sum_{i=2}^{k-1}\sum_{{|A|=|B|=k}\atop{|A\cap B|=i}}\Big(\mathbb{E}[X_AX_B]-\mathbb{E}[X_A]\mathbb{E}[X_B]\Big)\nonumber\\
&\le& \sum_{|A|=k} \mathbb{E}[X_A]+
\sum_{i=2}^{k-1}\sum_{{|A|=|B|=k}\atop{|A\cap B|=i}}\Big(\mathbb{P}[X_A=1]\mathbb{P}[X_B=1|X_A=1]-\mathbb{P}[X_A=1]\mathbb{P}[X_B=1]\Big)\nonumber\\
&\le& \mathbb{E}[X]+{n\choose k}p_k\cdot \sum_{i=2}^{k-1} {k\choose i}{{n-k}\choose{k-i}}(p_{k,i}-p_k)\label{l-cov}\ .
\end{eqnarray}

As a warm-up, we first show that a rather crude estimate for (\ref{l-cov}) suffices to prove that $\mathbb{V}ar[X]=o(\mathbb{E}^2[X])$ for $k=o(\sqrt{n})$. We start with the following bound for $p_{k,i}$.

\begin{lemma}\lab{claim3}
For $2\le i\le k-1$,
$$
p_{k,i} = O\left({{k-i}\choose{\frac{k-i}{2}}}^i2^{-(k-i)i}p_{k-i}\right)\ .
$$
Also,
$\frac{p_{k,i}}{p_k}\le Ce^{k\log\frac{k}{k-i}}$,
for a sufficiently large constant $C>0$.
\end{lemma}

\proof First, given $H$, expose the edges from $H$ to the
remaining $k-i$ vertices (denote the latter set by $X$). For every
$v\in H$, we require $d(v,X)=\kh-d_H(v)$. This happens with
probability
$$
{{k-i}\choose{\kh-d_H(v)}}2^{-k+i} \le {{k-i}\choose
{\frac{k-i}{2}}}2^{-k+i}
$$
(the middle binomial coefficient is the largest one). Hence the probability that all $i$ vertices from $V(H)$ have the required degree of $\kh$ in $G$ is at most the $i$-th power of the right hand side of the above expression.

Now, conditioned on the edges from $H$ to $X$, we ask what is the
probability that the subgraph spanned by $X$  has the required
degree sequence (each $v\in X$ should have exactly $\kh-d(v,H)$
neighbors in $X$). Observe that
by Lemma \ref{newclaim1}
 the probability that
$G[X]$ has the required degree sequence is at most $C_0p_{k-i}$ for
some absolute constant $C_0>0$, providing the first claimed estimate for $p_{k,i}$.

From Theorem~\ref{t:counting}, $p_t=\Theta\left(2^{-2{t\choose 2}}{{t-1}\choose{\lfloor (t-1)/2\rfloor}}^t\right)$. Therefore, the ratio
$p_{k,i}/p_k$ can be estimated as follows:
\begin{eqnarray*}
\frac{p_{k,i}}{p_k}&\le&  C_0{{k-i}\choose{\frac{k-i}{2}}}^i2^{-(k-i)i} \frac{p_{k-i}}{p_k}\le
C\left[\frac{{{k-i}\choose{\frac{k-i}{2}}}2^{-(k-i)}}
               {{{k-1}\choose{\frac{k-1}{2}}}2^{-(k-1)}}
   \right]^k \\
   &=&
   C\left[\frac{(\frac{k-1}{2})^2(\frac{k-1}{2}-1)^2\cdots (\frac{k-i}{2}+1)^2}
   {(k-1)(k-2)\cdots(k-i+1)}2^{i-1}\right]^k\\
   &=& C \left[ \frac{k-1}{k-2}\frac{k-3}{k-4}\cdots\frac{k-i+2}{k-i+1}\right]^k\\
   &\le& C\exp\left\{\Big(\frac{1}{k-2}+\frac{1}{k-4}+\cdots+\frac{1}{k-i+1}\Big)k\right\}\ .
\end{eqnarray*}
Observe that $\sum_{j=k-i+1}^{k-2}\frac{1}{j}< \int_{k-i}^{k}\frac{dx}{x}=\log\frac{k}{k-i}$. This completes the proof of the second part of the lemma.
\qed

\medskip

Now we complete a proof of a weaker version of the lower bound of Theorem \ref{th1}, by showing that
\whp\ $G(n,1/2)$ contains an induced $\kh$-regular subgraph on $k=o(\sqrt{n})$ vertices.
Omitting the negative term of $p_k$ in the sum in (\ref{l-cov}) and using $\mathbb{E}[X]={n\choose k}p_k$, we obtain:
\begin{equation}\label{VE}
\frac{\mathbb{V}ar[X]}{\mathbb{E}^2[X]} \le \frac{\sum_{i=2}^{k-1}{k\choose i}{{n-k}\choose{k-i}}p_{k,i}}{{n\choose k}p_k} +\frac{1}{\mathbb{E}[X]} =
\sum_{i=2}^{k-1}\frac{{k\choose i}{{n-k}\choose{k-i}}}{{n\choose k}} \frac{p_{k,i}}{p_k}+o(1) .
\end{equation}
Denote
$$
g(i)=\frac{{k\choose i}{{n-k}\choose{k-i}}}{{n\choose k}} \frac{p_{k,i}}{p_k}\,.
$$
Let us first estimate the ratio of the binomial coefficients involved in the definition of $g(i)$.
\begin{eqnarray*}
\frac{{k\choose i}{{n-k}\choose{k-i}}}{{n\choose k}}&\le& \frac{\left(\frac{ek}{i}\right)^i{n\choose{k-i}}}{{n\choose k}}
= \left(\frac{ek}{i}\right)^i\, \frac{k(k-1)\cdots (k-i+1)}{(n-k+i)(n-k+i-1)\cdots(n-k+1)}\\
&\le& \left(\frac{ek}{i}\right)^i\, \left(\frac{k}{n-k+i}\right)^i
\le \left(\frac{3k^2}{in}\right)^i\ .
\end{eqnarray*}

To analyze the asymptotic behavior of $g(i)$, we consider three cases.

{\bf Case 1. $i\le k/2$.} In this case, by Lemma \ref{claim3} and the inequality $\log(1+x)\le x$ for $x\ge 0$ we have:
$$
\frac{p_{k,i}}{p_k}\le Ce^{k\log\frac{k}{k-i}}= Ce^{k\log(1+\frac{i}{k-i})}\le Ce^{k\frac{i}{k-i}}
\le Ce^{2i}\,.
$$
We thus get the following estimate for $g(i)$:
$$
g(i)=\frac{{k\choose i}{{n-k}\choose{k-i}}}{{n\choose k}} \frac{p_{k,i}}{p_k}\le
\left(\frac{3k^2}{in}\right)^i\cdot Ce^{2i} \le C\left(\frac{3e^2k^2}{in}\right)^i= (o(1))^i\ .
$$

{\bf Case 2. $k/2\le i\le k-\frac{k}{\log k}$.}
Recalling Lemma \ref{claim3} again, we have
${p_{k,i}}/{p_k}\le Ce^{k\log\frac{k}{k-i}}\le Ce^{k\log\log k}$.  Hence in this case
\begin{eqnarray*}
g(i) &\le& \left(\frac{3k^2}{in}\right)^i\, Ce^{k\log\log k}\le
\left(\frac{6k}{n}\right)^{k/2}\, Ce^{k\log\log k}\\
&\le& \left(\frac{1}{\sqrt{k}}\right)^{k/2}\, Ce^{k\log\log k}
= Ce^{-\frac{k\log k}{4}+k\log\log k} \le e^{-k}\,.
\end{eqnarray*}
For future reference, it is important to note here that in the calculation above we used $6k/n\le k^{-1/2}$. This inequality stays valid as long as $k\le (n/6)^{2/3}$.

{\bf Case 3. $i\ge k-\frac{k}{\log k}$.}
In this case it suffices to use the trivial estimate $p_{k,i}\le 1$. We also need that $\mathbb{E}[X]={n\choose k}p_k=(\omega(1))^k$.
Therefore,
$$
g(i)= \frac{{k\choose i}{{n-k}\choose{k-i}}p_{k,i}}{{n\choose k}p_k} \le \frac{{k\choose i}{{n-k}\choose{k-i}}}{(\omega(1))^k}\le
\frac{2^kn^{k-i}}{(\omega(1))^k}\le \frac{2^kn^{k/\log k}}{(\omega(1))^k}=\frac{e^{O(k)}}{(\omega(1))^k}\le e^{-k}\,.
$$
In the above calculation we used the assumption $\log n=O(\log k)$. In the complementary case $k=n^{o(1)}$ the expression ${n\choose k}p_k$ behaves like $\big(cn/k^{3/2}\big)^k\ge n^{k/2}$, while the numerator in the expression for $g(i)$ is at most $2^kn^{k/\log k}=n^{o(k)}$, and the estimate works as well.

\medskip

Now we proceed to the proof of the ``real" lower bound of Theorem \ref{th1}, i.e. assume that $k$ satisfies $k=o(n^{2/3})$. In this case estimating the variance of the random variable $X$, defined as the number of induced $\kh$-regular subgraphs on $k$ vertices, becomes much more delicate. We can no longer ignore the negative term of $p_k$ in the sum in (\ref{l-cov}). Instead, we show that for small values of $i$ in this sum $p_{k,i}$ is asymptotically equal to $p_k$. In words, this means that knowing the edges spanned by the first $i$ vertices of a random graph $G=G(k,1/2)$  does not affect by much the probability of $G$ being $\kh$-regular. We claim this formally for $i=o\sqrt{k}$ in the following key lemma.

\begin{lemma}\label{conj2}
For $i=o(\sqrt{k})$,
$$
p_{k,i}=(1+o(1))p_k .
$$
\end{lemma}
The proof of this lemma is rather involved technically. We thus postpone it to the next section. We now show how to complete the proof assuming its correctness. We first repeat estimate (\ref{l-cov}):
\begin{eqnarray}
\frac{\mathbb{V}ar[X]}{\mathbb{E}^2[X]} &\le& \frac{1}{\mathbb{E}[X]}+
\sum_{i=2}^{k-1}\frac{{k\choose i}{{n-k}\choose{k-i}}(p_{k,i}-p_k)}{{n\choose k}p_k}\nonumber\\
&\le&
o(1)+\sum_{i=2}^{t} \frac{{k\choose i}{{n-k}\choose{k-i}}(p_{k,i}-p_k)}{{n\choose k}p_k}
+ \sum_{i=t}^{k-1}\frac{{k\choose i}{{n-k}\choose{k-i}}}{{n\choose k}} \frac{p_{k,i}}{p_k}\,,\label{lb10}
\end{eqnarray}
where $t=t(k,n)$ is chosen so that $t=\omega(k^2/n)$ but $t=o(\sqrt{k})$. Since $k=o(n^{2/3})$ such a function is easily seen to exist. Due to our choice of $t$ we can apply Lemma \ref{conj2} to the first sum above. It thus follows that
$$
\sum_{i=2}^{t}\frac{{k\choose i}{{n-k}\choose{k-i}}(p_{k,i}-p_k)}{{n\choose k}p_k} =
\sum_{i=2}^{t}\frac{{k\choose i}{{n-k}\choose{k-i}}\cdot o(p_k)}{{n\choose k}p_k}
\le o(1)\cdot \frac{\sum_{i=0}^{k}{k\choose i}{{n-k}\choose{k-i}}}{{n\choose k}} = o(1)\,.
$$
As for the second sum in (\ref{lb10}) we can utilize the same case analysis as done before for $k=o(\sqrt{n})$. The only difference is in Case 1, that now covers all $i$ from $t$ till $k/2$. Therefore, for every $i$ in this new interval we have
$$
\frac{{k\choose i}{{n-k}\choose{k-i}}}{{n\choose k}} \frac{p_{k,i}}{p_k}
\le C\left(\frac{3e^2k^2}{in}\right)^i\le C\left(\frac{3e^2k^2}{tn}\right)^i= (o(1))^i\ .
$$
This completes the proof of Theorem \ref{th1}.
\qed

\section{Proof of key lemma}

The proof of Lemma \ref{conj2} is overall along the lines of the proof of Lemma~\ref{claim3}, though
requiring a much more detailed examination of the  probabilities involved. Let $\kk$ be odd and,
for simplicity, denote $(\kk-1)/2$ by $d$. Let $\Dscr_i$ be the set of  integer vectors
$\dv=(d_1,\ldots, d_i)$ such that $0\le d_j\le i-1$ for $1\le j\le i$, and $\sum_j d_j$ is even.
Given $\dv\in \Dscr_i$, let   $N(\dv)$ denote the
 number of  graphs $G$ on vertex set $[\kk]$ for which  $G[i]$ has no edges, and $d_G(j)=d-d_j$
for $j\in [i]$. Note that if $\dv$ is the degree sequence of a graph $H$ on vertex set $[i]$, then
$N(\dv)$ is the number of $d$-regular graphs $G$ on vertex set $[\kk]$ for which  $G[i]=H$.

\begin{prop}\lab{p:ratio}
Assume $i=o(\sqrt \kk)$. Given $\dv\in \Dscr_i$ and a nonnegative vector
$\sv=(s_1,\ldots, s_i)$, put $\dv'=\dv-\sv$.
Then, uniformly over such
$\dv$ and  $\sv$ with the additional properties that  $\dv'\in \Dscr_i$ and $\sum_{j=1}^i s_j\le k^{3/4}$,
$$
\frac{N(\dv)}{N(\dv')}\sim
\frac{\displaystyle   \prod_{j=1}^i {d-d_j+s_j \choose s_j}}
{\displaystyle   \prod_{j=1}^i{d+1-(i-d_j) \choose s_j}}.
$$
\end{prop}

\proof 
We use a comparison type argument. Since it is quite complicated, we give the idea of the proof first.  
For any vector $\cv=(c_1,\ldots, c_j)$, write $\cv^*$ for the vector $(d-c_1,\ldots, d-c_j)$.
Let $V_1=\{1,\ldots, i\}$ 
and $V_2=\{i+1,\ldots, \kk\}$. For simplicity, suppose that $s_1=s_2=1$, and $s_j=0$ for $j\ge 3$.  We can compute $N(\dv)$ as the number 
of possible outcomes of two steps. The first step is to choose a bipartite graph $B$ with bipartition $(V_1,V_2)$ and degree sequence 
$\dv^*$ in $V_1$. The second step is to add the remaining edges between vertices in $V_2$ such that those vertices will have degree $d$. 
By comparison, to count $N(\dv')$ we choose in the first step $B'$ with degree sequence 
$\dv'^*$ in $V_1$, and then do the second step for 
each 
such $B'$. The proof hinges around the fact that there is a correspondence between the set of possible $B$ and $B'$ such that the number 
of ways of performing the second step is roughly the same, at least for most of the corresponding pairs $(B,B')$.

The correspondence is many-to-many. For a graph $B$ we may add two edges, incident with vertices 1 and 2, to obtain a graph $B'$.  The 
number of ways this can be done, without creating multiple edges, is $ 
\prod_{j=1}^2\big(k-i-(d-d_j)\big)=\prod_{j=1}^2\big(d+1-(i-d_j)\big)$.  Conversely, each $B'$ comes from $\prod_{j=1}^2 (d -d_j+1)\big)$ 
different $B$. The ratio of these quantities gives the asymptotic ratio between $N(\dv)$ and $N(\dv')$ claimed in the theorem. Our actual 
argument gets more complicated because not only some bipartite graphs must be excluded, but also some sets of edges to be added to them. 
So we will present equations relating to the above argument in a slightly different form to make exclusion of various terms easier.

Let $\Bscr$ denote the set of bipartite graphs   with bipartition
$(V_1,V_2)$.  For $B\in \Bscr$, write
${\bf D}_j(B)$ for the degree sequence of $B$ on the vertices in $V_j$ (in non-decreasing order), so
${\bf D}_1(B)=\big(d_B(1),\ldots, d_B(i)\big)$ and ${\bf D}_2(B)=\big(d_B(i+1),\ldots,
d_B(\kk)\big)$. For $\dv=(d_1,\ldots, d_i)$,  let
$\Bscr(\dv^*)$ denote $\{ B\in\Bscr : {\bf D}_1(B)=\dv^*\}$.  
Let $G\big({\bf D}_2(B)^*\big)$ denote the number of graphs with degree sequence 
${\bf D}_2(B)^*$. Clearly, if $\dv$ is graphical,
\bel{Nd}
N(\dv) = \sum_{B\in \Bscr(\dv^*)}G\big({\bf D}_2(B)^*\big).
\ee

Suppose that we wish to add to $B$ a set $S$ of   edges joining $V_1$ and $V_2$,  without creating any
multiple edges,  such that the degree of $j\in V_1$ ($1\le j\le i$) in the graph induced by $S$ is
$s_j$ (as given in the statement of the proposition). The family of all such sets
$S$ will be denoted by $\Sscr(B,\sv)$. Note that necessarily   $|S|\le k^{3/4}$ for $S\in \Sscr(B,\sv)$.  The cardinality of $\Sscr(B,\sv)$ is $\prod_{j=1}^i{d+1-(i-d_j)
\choose s_j}$, because $d_B(j)=d_j^* =d-d_j$, so (as in the sketch above) $j$ has $d+1-(i-d_j)$ spare vertices in 
$V_2$
to which it may be joined. Hence we can somewhat artificially rewrite~\eqn{Nd} as
\bel{Nd2}
N(\dv) = \frac{1}{\displaystyle\prod_{j=1}^i{d+1-(i-d_j) \choose s_j}} \sum_{B\in
\Bscr(\dv^*)}\sum_{S\in\Sscr(B,\sv)} G\big({\bf D}_2(B)^*\big).
\ee
Also for ${B'\in \Bscr(\dv'^*)}$
define $\Sscr'(B',\sv)$ to be the family of sets $S\subseteq E(B')$ such that   the degree of
$j\in V_1$ in the graph induced by $S$ is
$s_j$  ($1\le j\le i$).  Since $d_{B'}(j)=d-d_j'$ and $d_j'=d_j-s_j$, a similar argument gives
\bel{Ndp}
N(\dv') = \frac{1}{\displaystyle\prod_{j=1}^i {d-d_j+s_j \choose s_j}} \sum_{B'\in
\Bscr(\dv'^*)}\sum_{S\in\Sscr'(B',\sv)} G\big({\bf D}_2(B')^*\big).
\ee
The rest of the proof consists of showing that the significant terms in the last two
equations can be put into 1-1 correspondence such that corresponding terms are asymptotically
equal.

We first need to show that for a typical $B\in\Bscr(\dv^*)$, the variance
of the elements of ${\bf D}_2(B)$ (as a sequence) is small.
Given $\dv$,
define $\dbar=(\kk-i)^{-1}\sum_{j\in V_1} (d-d_j)$, and note that this
is equal to
$(\kk-i)^{-1}\sum_{j\in V_2} d_B(j)$ for every $B\in \Bscr $.
Since $\dbar$ is determined
uniquely by $\dv$, the value of $\dbar$ is the same
for all $B\in \Bscr(\dv^*)$.
\begin{lemma}\lab{l:seqvar}
Let $\dv\in\Dscr_i$, and select $B$ uniformly at random from $\Bscr(\dv^*)$.
Then
$$
\ex\left( \sum_{j\in V_2} \big(\dbar   -d_B(j)\big)^2\right) \le i(k-i).
$$
\end{lemma}
\proof
First observe that in $B$, the neighbours of any vertex $t \in V_1$ form a
random subset of $V_2$ of size $d_t^*$, and these subsets are independent for different $t$. So
for fixed $j\in V_2$, $d_B(j)$ is distributed as a sum of $i$ independent 0-1 variables with
mean $\sum_{t\in V_1} d_t^*/(k-i)=\dbar$. It follows that  the  variance of $d_B(j)$ 
is less than $i$. 
Hence  $\ex \big(\dbar   -d_B(j)\big)^2 <i$, and the lemma follows by linearity of
expectation.
\qed
\smallskip

Returning to the proof of the proposition,
we will apply Theorem~\ref{t:counting} to estimate $G\big({\bf D}_2(B)^*\big)$. 
This graph has $k-i$ vertices, degree sequence $\{d-d_B(j), j \in V_2\}$, and its  
number of edges is $e\big(G({\bf D}_2(B)^*)\big)=(k-i)d-e(B)$, where
$e(B)=\sum_{j \in V_2} d_B(j)=(k-i)\dbar$ is the number of edges in bipartite graph $B$.
Consider $\la$ from Theorem~\ref{t:counting}. Using the representation of $\dbar$ as the average degree 
of $B$ in $V_2$, we see that
\bel{lambda}
\la = \la(\dv):=\frac{e\big(G({\bf D}_2(B)^*)\big)}{(k-i)(k-i-1)}=\frac{d-\dbar}{k-i-1}.
\ee
The product of binomials in~\eqn{formula} is in this case
\bel{binprod}
\prod_{j\in V_2}  {k-i-1\choose d-d_B(j)}.
\ee
For every $B\in \Bscr(\dv^*)$, all components of the vector $ {\bf D}_2(B)$ are at most $|V_1|=i$. 
Thus 
$$x_j:=\frac{k-i-1}{2} - (d-d_B(j))=d_B(j)-i/2=O(i)=o(\sqrt k).$$ We have
\bel{binapprox}
{a \choose a/2+x}={a \choose \lfloor a/2\rfloor}\exp\big(-2x^2/a+O(x^3/a^2)\big)
\ee
for $x=o(\sqrt a)$, which may be established 
for instance by analyzing the ratio of the binomial coefficients. Hence
\begin{equation}
\label{binomials}
\prod_{j=i+1}^{k} {k-i-1\choose d-d_B(j)} =
{{k-i-1}\choose{\frac{k-i-1}{2}}}^{k-i}
\exp\left( \frac{-2\sum x_j^2}{k-i}+o(i)\right).
\end{equation}
(Note that here and in the rest of the proof, the asymptotic relations hold uniformly over $\dv\in \Dscr_i$.)
Since $i=o(\sqrt{k})$ we can choose a function $\omega$ of $n$ such that $\omega\to\infty$ and $ 
\omega^2 i = o(\sqrt k)$. Define 
$\hat \Bscr_\omega(\dv^*)$  to be the subset of  $\Bscr(\dv^*)$ that contains those $B$ for which
\bel{hatdef}
\sum_{j\in V_2} \big(\dbar   -d_B(j)\big)^2 \le \omega i(\kk-i).
\ee
Since $\sum_{j \in V_2}d_B(j)=e(B)$ is the same for all bipartite graphs $B \in \Bscr(\dv^*)$, by 
definition of $x_j$ we 
have that $\sum_j x_j^2 - \sum_j d^2_B(j)$ also does not depend on $B$. Similarly, the sum in~\eqn{hatdef} differs from  $\sum_j d^2_B(j)$ by a constant independent of $B$.
Therefore $\sum_j x_j^2$ for all $B\in \hat \Bscr_{\omega/2}(\dv^*)$ is smaller than the 
corresponding sum for $B\in \Bscr(\dv^*)\setminus \hat \Bscr_\omega(\dv^*)$ by an additive term
of at least $\omega i(\kk-i)/2$. This implies that 
the product of binomials in (\ref{binomials}) is larger,
for all $B\in \hat \Bscr_{\omega/2}(\dv^*)$,
than for any  $B\in \Bscr(\dv^*)\setminus \hat \Bscr_\omega(\dv^*)$.
Also, from   Lemma~\ref{l:seqvar} and Markov's inequality, almost all members of  $\Bscr(\dv^*)$ are in  
$\hat \Bscr_{\omega/2}(\dv^*)$. Moreover, since all degrees in degree sequence ${\bf D}_2(B)^*$ deviate 
from $(k-i)/2$ by at most $O(i)=o(\sqrt{k})$, the function
$f({\bf D}_2(B)^*)$ from Theorem~\ref{t:counting} is $\, \sim\, \sqrt{2}e^{1/4}$ for all $B \in 
\Bscr(\dv^*)$.
Combining these observations, we conclude that the contribution to~\eqn{Nd} from 
$B\notin  \hat \Bscr_\omega(\dv^*)$ is $o\big(N(\dv)\big)$. 
Thus, the same observation holds for~\eqn{Nd2}. That is,
\bel{Nd3}
N(\dv) \sim \frac{1}{\displaystyle\prod_{j=1}^i{d+1-(i-d_j) \choose s_j}} \sum_{B\in
\hat \Bscr_\omega(\dv^*)}\sum_{S\in\Sscr(B,\sv)} G\big({\bf D}_2(B)^*\big).
\ee
We also note for later use, that by~\eqn{hatdef} and Cauchy's inequality, for all $B\in \hat 
\Bscr_\omega(\dv^*)$
\bel{squares}
\sum_{j\in V_2}  \big|\dbar   -d_B(j)\big| \le
(k-i)\sqrt{\omega i}.
\ee

Fix $B\in \hat \Bscr_\omega(\dv^*)$.
Consider $S$ chosen uniformly at random from $\Sscr(B,\sv)$, and let $r_m(S)$ 
denote the number of edges of $S$ incident with a vertex $m\in V_2$.  Fixing $m$
and using that $|S| \leq k^{3/4}$, we can bound
the probability that $r_m(S)\ge 5$ by 
$$\sum_{\{j_1,\ldots , j_5\}\subseteq V_1}\, \prod_{t=1}^5 \frac{s_{j_t}}{k-i-(d-d_{j_t})}
\le  \left(\sum_{j\in V_1} \frac{s_{j}}{d-i}\right)^5
= \left(\frac{|S|}{d-i}\right)^5 = O(k^{-5/4}).$$
Hence by Markov's inequality, with probability $1-O(k^{-1/4})$, $S\in \Sscr(B,\sv)$ satisfies
\medskip

\noindent
(i) $\max_{j\in V_2} r_j(S) \le 4$.
\medskip

Note also that each vertex of $V_1$ has, as crude bounds, between $(k-i)/3$ and $2(k-i)/3$ vertices of 
$V_2$ eligible to choose for an edge of $S$ (at least, for large $k$).  Hence, the expected value of 
$\big|\dbar -d_B(j)\big|$ amongst all such vertices is at most $3\sqrt{\omega i}$ by~\eqn{squares}. Note 
that we may choose the edges in $S$ incident with any given vertex sequentially, each time selecting a 
random neighbour from those vertices of $V_2$ still eligible to be joined to. For each such edge joining 
to such a random vertex $j\in V_2$, the unconditional expected value of $\big|\dbar -d_B(j)\big|$ is at 
most $3\sqrt{\omega i}$. Thus by Markov's inequality, and noting that $\sum r_j(S) \le k^{3/4}$, we 
deduce that almost all $S\in \Sscr(B,\sv)$ (more precisely the fraction
$1-3/\sqrt{\omega}=1-o(1)$  of them, at least) satisfy \medskip

\noindent
(ii) $\displaystyle
\sum_{j\in V_2}  \big|\dbar   -d_B(j)\big|  r_j(S) \le  \omega \sqrt i k^{3/4}.
$
\medskip

Define $\hat\Sscr(B,\sv)$ to be the set of   $S\in \Sscr(B,\sv)$ satisfying both the properties  (i) and (ii). Then, since each $S\in \Sscr(B,\sv)$ contributes equally to~\eqn{Nd3},
\bel{Nd4}
N(\dv) \sim \frac{1}{\displaystyle\prod_{j=1}^i{d+1-(i-d_j) \choose s_j}} \sum_{B\in
\hat \Bscr_\omega(\dv^*)}\sum_{S\in\hat\Sscr(B,\sv)} G\big({\bf D}_2(B)^*\big).
\ee

Let $B\in \hat \Bscr_\omega(\dv^*)$ and $S\in\hat\Sscr(B,\sv)$. Then, using (\ref{hatdef}) together with (i) and (ii), we get 
\bean
\sum_{j\in V_2} \big(\dbar   -d_B(j)-r_j(S)\big)^2
&\le&
\omega i(k-i) +
2 \sum_{j\in V_2}  \big|\dbar   -d_B(j)\big|  r_j(S) +
\sum_{j\in V_2} r_j(S)^2\\
&\leq& \omega i(k-i)+2\omega\sqrt{i}k^{3/4}+16(k-i)    \sim
\omega i\kk.
\eean
Hence, for $k$ sufficiently large, those $S$ appearing in the range of the summation in~\eqn{Nd4} satisfy $B+S\in
\hat \Bscr_{2\omega} (\dv'^*)$, where $B+S$ is the graph obtained by adding the edges in $S$ to $B$ (and noting that $\dv'^*=\dv^* +\sv$). Since, as we saw,   the contribution to~\eqn{Nd} from $B\notin  \hat \Bscr_\omega(\dv^*)$ is $o\big(N(\dv)\big)$, we may also relax the constraint on $B$ in the summation in~\eqn{Nd4}, to become $B \in \hat \Bscr_{2\omega}(\dv^*)$. Now redefining $2\omega$ as $\omega$, we obtain
\bel{Nd5}
N(\dv) \sim \frac{1}{\displaystyle\prod_{j=1}^i{d+1-(i-d_j) \choose s_j}} \sum_{(B,S)\in \Wscr} G\big({\bf D}_2(B)^*\big),
\ee
where $\Wscr$ denotes the set of all $(B,S)$ such that
$B\in \hat \Bscr_\omega(\dv^*)$, $S\in\hat\Sscr(B,\sv)$ and $B+ S\in \hat \Bscr_\omega(\dv'^*)$.

Define $\hat\Sscr'(B',\sv)$, analogous to $\hat\Sscr(B,\sv)$, to be the set of $S\in \Sscr'(B',\sv)$ with maximum degree in $V_2$ at most 5 and also obeying property (ii) above, where $B=B'-S$. Then the above argument applied to~\eqn{Ndp}, with suitable small modification, gives
\bel{Ndp5}
N(\dv') \sim \frac{1}{\displaystyle\prod_{j=1}^i {d-d_j+s_j \choose s_j}}
\sum_{(B',S)\in \Wscr'}  G\big({\bf D}_2(B')^*\big)
\ee
where $\Wscr'$ denotes the set of all $(B',S)$ such that
$B'\in \Bscr(\dv'^*)$, $S\in\hat\Sscr(B',\sv)$ and $B'- S\in \hat \Bscr_\omega (\dv ^*)$.

Observe that $(B,S)\in \Wscr$ if and only if $(B+S,S)\in \Wscr'$. So the summation in~\eqn{Ndp5} is equal to
$$
\sum_{(B,S)\in \Wscr} G\big({\bf D}_2(B+S)^*\big).
$$
Hence, comparing with~\eqn{Nd5}, the proposition follows if we show that
\bel{toshow}
G\big({\bf D}_2(B)^*\big)\sim G\big({\bf D}_2(B+S)^*\big)
\ee
uniformly for all $(B,S)\in \Wscr$.

We may apply~\eqn{formula} to both sides of~\eqn{toshow}. Write
$g(\la,n)= \big(\la^{\la} (1-\la)^{1-\la}\big)^{n\choose 2}.$ Notice that $\la(\dv)$, as defined in \eqn{lambda}, satisfies:
$$  
\la(\dv)=\frac{\sum_{j\in V_2}(d-d_B(j))}{(k-i)(k-i-1)}
$$
--- which is exactly $\la$ for the degree sequence ${\bf D}_2(B)^*$ as defined in Theorem \ref{t:counting}. 
The same applies to $\la(\dv')$ and the degree sequence ${\bf D}_2(B+S)^*$. Using
$$
\sum_{j\in V_2}(d-d_B(j))=(k-i)d-\sum_{t \in V_1}d_B(t)=(k-i)d-id+\sum_{t=1}^i 
d_t=\frac{(k-2i)(k-1)}{2}+\Theta(i^2),
$$
it is easy to derive from \eqn{lambda} that both $\la(\dv)$ and $\la(\dv')$ are 
$1/2+O(i^2/k^2)=1/2+o(k^{-1})$ for all $\dv,\dv'\in \Dscr_i$. 
For such $\la$ the derivative of $\log g(\la,k-i)$ is $\displaystyle{k \choose 2}\log (\la/(1-\la))=O(k^2i^2/k^2)=o(k)$.   Moreover, $|\la(\dv)-\la(\dv')| = O(k^{-5/4})$ 
because the values of $\dbar$ for $B$ and $B+S$ differ by 
$O(|S|/k)$. Hence  $g(\la(\dv),k-i)\sim g(\la(\dv'),k-i)$. It is now also easy to see that 
$f({\bf D}_2(B)^*), f({\bf D}_2(B+S)^*)$ from Theorem \ref{t:counting} satisfy: 
$f({\bf D}_2(B)^*)\sim f({\bf D}_2(B+S)^*)\sim \sqrt{2}e^{1/4}$, since all degrees in these two degree sequences deviate from 
$(k-i)/2$ by $O(i)=o(\sqrt{k-i})$.

It only remains to consider the product of binomials in the two sides of~\eqn{toshow}. 
Recalling the expression~\eqn{binprod},  the ratio of these two products in the two cases is
$$
\prod_{j \in V_2} {k-i-1\choose d-d_B(j)}\Big/{k-i-1\choose d-d_B(j)-r_j(S)}.
$$
Since
$B\in \hat \Bscr_\omega(\dv^*)$, all $r_j(S) \leq 4$  and $\sum_j r_j(S) \leq k^{3/4}$, so using 
$d=(k-1)/2$ and $d_B(j)\le i$, 
this expression is, up to a multiplicative factor of 
$1+O(k^{-1}) \sum_j r_j^2(S) =1+o(1)$, equal to
\bean
\prod_{j \in V_2} \left(\frac{k-i-1-d+d_B(j)}{d-d_B(j)}\right)^{r_j(S)}
&=&
\prod_{j \in V_2} \left(\frac{d-i/2 -\big(i/2 -d_B(j)\big)}{d-i/2+\big(i/2-d_B(j)\big)}\right)^{r_j(S)}\\
&=& \prod_{j \in V_2}  \left(1+O\left(\frac{i/2-d_B(j)}{k}\right)\right)^{r_j(S)}.
\eean
By its definition, $\dbar=(k-i)^{-1}\sum_{j \in V_1}(d-d_j)=i/2+O(i^2/k)=i/2+o(1)$, and so using 
condition (ii) (the right hand side of which is $\omega\sqrt{i}k^{3/4}=o(k)$), 
and not forgetting $\sum r_j(S)\le k^{3/4}$, we get
$$\sum_{j \in V_2} r_j(S) \frac{|i/2-d_B(j)|}{k} \leq 
\sum_{j \in V_2}\frac{|\dbar-d_B(j)|r_j(S)}{k}+\sum_{j \in V_2}\frac{|\dbar-i/2|r_j(S)}{k} =o(1).$$
Hence, the expression above is asymptotic to 1. This argument shows that~\eqn{toshow} holds with the 
required uniformity.
\qed

For a slightly simpler version of the formula in Proposition~\ref{p:ratio}, put
$$
\dhat = d-\frac{1}{2}(i-1) = \frac12(k-i)
$$
(which is in some sense the average degree of vertices of side $V_1$ in the bipartite graph $B$) and
$$
\delta_j = d_j-\frac{1}{2}(i-1).
$$
Then the proposition gives
$$\frac{N(\dv)}{N(\dv')}\sim
\displaystyle   \prod_{j=1}^i
\frac{[\dhat-\delta_j+s_j ]_{s_j}}{[\dhat+\delta_j]_{s_j}}.
$$
Recalling that  $\delta_j$ and $s_j$ are at most
$i=o(\sqrt k)$ and using $\log(1+x)=x-x^2/2+O(x^3)$, we have
\begin{eqnarray*}
[\dhat-\delta_j+s_j ]_{s_j}
&=&
\dhat^{s_j}\,\exp\big( - \delta_js_j/\dhat +s_j^2/2\dhat+o(1/\sqrt k)\big),\\
{}[\dhat+\delta_j]_{s_j}
&=&\dhat^{s_j}\,\exp\big(   \delta_js_j/\dhat -s_j^2/2\dhat+o(1/\sqrt k)\big).
\eean
Thus, we may rewrite the assertion of Proposition~\ref{p:ratio} as
\bel{simpler}
\frac{N(\dv)}{N(\dv')}\sim
\exp\left\{\sum_{j=1}^i \big(-2\delta_js_j/\dhat +s_j^2/\dhat\,\big)\right\}.
\ee
To proceed, we extend this formula  so that $\sv$ is permitted to have negative entries.
\begin{corollary}\lab{c:expratio}
Assume $i=o(\sqrt \kk)$. Given $\dv\in \Dscr_i$ and an integer vector
$\sv=(s_1,\ldots, s_i)$, put $\dv'=\dv-\sv$.
Then, uniformly over such
$\dv$ and  $\sv$ with the additional properties that  $\dv'\in \Dscr_i$ and $\sum_{j=1}^i |s_j|\le k^{3/4}$,
$$
\frac{N(\dv)}{N(\dv')}\sim
\exp\left\{  \frac{1}{\dhat}\sum_{j=1}^i \big(-2\delta_js_j +s_j^2\big)\right\}.
$$
\end{corollary}
\proof
Define the vector $\sv'$ by turning the negative entries of $\sv$ into 0; 
that is, the   $j$th entry of $\sv'$ is $s_j$ if $s_j\ge 0$, and 0 otherwise. 
Let $\sv''=\sv'-\sv$. The $j$th entry of $\sv''$ is $-s_j$ if $s_j< 0$, and 0 otherwise.
We can now estimate the product
$$
\frac{N(\dv)}{N(\dv-\sv')}\cdot \frac{N(\dv - \sv')}{N(\dv')}=
\frac{N(\dv)}{N(\dv-\sv')} \Big/ \frac{N(\dv')}{N(\dv'-\sv'')}
$$
using two applications of (\ref{simpler}). 
First,  
$$\frac{N(\dv)}{N(\dv-\sv')} \sim 
\exp\left\{  \frac{1}{\dhat}\sum_{j=1}^i \big(-2\delta_js'_j +(s'_j)^2\big)\right\}=
\exp\left\{ \frac{1}{\dhat}\sum_{s_j \geq 0} \big(-2\delta_js_j 
+s_j^2\big)\right\}.$$
Next, note that all entries of $\sv''=\sv'-\sv$ are nonnegative and that
the $\delta'_j$ defined for degree sequence $\dv'$ equals $\delta_j-s_j$.
Applying (\ref{simpler}) again,
$$\frac{N(\dv')}{N(\dv'-\sv'')} \sim
\exp\left\{  \frac{1}{\dhat}\sum_{j=1}^i \big(-2\delta'_js''_j +(s''_j)^2\big)\right\}=
\exp\left\{ \frac{1}{\dhat}\sum_{s_j < 0} \big(2(\delta_j-s_j)s_j 
+s_j^2\big)\right\}.$$
To complete the proof, divide the first formula by the second. \qed
\smallskip

Define $\dv_0$ to be the constant sequence of length $i$, 
all of whose entries are $\lfloor (i-1)/2\rfloor$.  
We can use the following result to compare the number of graphs with  an arbitrary 
degree sequence $\dv$ on $G[i]$ to the number  with $\dv_0$. Recall, however, 
that $N(\dv)$ is defined even if $\dv$ is not the degree sequence of any graph.

\begin{corollary}\lab{c:all}
\noindent (i) If $\dv\in \Dscr_i$ then $N(\dv)\le  N(\dv_0)\big(1+o(1)\big)$.

\noindent (ii) If, in addition, $\sum_{j=1}^i\delta_j^2=o(k)$, then $N(\dv)\sim N(\dv_0)$.
\end{corollary}
\proof The second part follows immediately from Corollary~\ref{c:expratio} by putting 
$s_j=\lceil\delta_j\rceil$ for each $j$, since if 
$\sum_{j=1}^i\delta_j^2=o(k)$ then by Cauchy's inequality $\sum |\delta_j|= o(k^{3/4})$. (Note 
also that $\dhat\sim k/2$.)

For the first part, let $\dv$ maximise $N(\dv)$. If $\sum_{j=1}^i\delta_j^2<k$ say, the above argument shows that $N(\dv)\le 
N(\dv_0)\big(1+o(1)\big)$. So assume that $\sum_{j=1}^i\delta_j^2>k$. Putting $s_j=\lceil \delta_j\rceil$ and applying 
Corollary~\ref{c:expratio} shows the result, provided $\sum |\lceil \delta_j\rceil|\le k^{3/4}$. If the latter condition fails, we can 
simply define $s_j=\alpha_j\delta_j$ for some $0 \leq \alpha_j \leq 1$ such that $\sum_j|s_j|$
is just below  $k^{3/4}$ and is even.
Since $|s_j| \leq |\delta_j|$ and they both have the same sign, we can conclude that $\sum_j (-2\delta_js_j+s_j^2 )\leq -\sum_js_j^2$.
By Cauchy's inequality, the sum of the squares of $s_j$ grows asymptotically faster than $k$.
Let $\sv=(s_1, \ldots, s_i)$ and let $\dv'=\dv-\sv$. Then, by Corollary \ref{c:expratio} 
we obtain $N(\dv)=o(N(\dv'))$, which contradicts the maximality assumption and proves the result. \qed

\smallskip

Define $\dv_1$ to be the constant sequence of length $k$, all of whose entries are $d=(k-1)/2$. We can 
now determine the asymptotic value of $N(\dv_0)$.
\begin{corollary}\lab{c:Nd0}
$ N(\dv_0)\sim G(\dv_1)/2^{i\choose 2}$.
\end{corollary}
\proof
Let $H$ be one of the $2^{i\choose 2}$ graphs on vertex set $[i]$ chosen at random, and let $\dv_H=
\{d_1, \ldots, d_i\}$ be 
its degree sequence. Then $d_j$ is a binomially distributed random variable with 
expectation $(i-1)/2$ and variance $(i-1)/4$. Hence for $\delta_j=d_j-(i-1)/2$ we have
$\mathbb{E}[\delta_j^2]=\mathbb{V}ar[d_j]=(i-1)/4$. Then
$$ 
\ex\sum_{j\in [i]} \delta_j ^2 \le i^2\,, $$ 
and, by Markov's inequality, \whp\, $\sum_{j\in [i]} \delta_j ^2=o(\kk)$. 
Thus, from Corollary~\ref{c:all}(ii) it follows that
for almost all graphs $H$, $N(\dv_H)\sim N(\dv_0)$. Part (i) of the 
same corollary shows that for all other graphs, $N(\dv_H)\le (1+o(1)) N(\dv_0)$. Since $N(\dv_H)$ is 
the number of $d$-regular graphs $G$ on vertex set $[\kk]$ for which 
$G[i]=H$, we have that
$G(\dv_1)=\sum_{H}N(\dv_H)=(1+o(1))N(\dv_0)2^{i\choose 2}$, and 
the corollary follows.  \qed

\medskip

\noindent
{\bf Proof of Lemma \ref{conj2}.}
From Corollary~\ref{c:all},
$$
p_{k,i}\sim  \pr[\mbox{$G(k,0.5)$ is $\kh$-regular} \mid
G[i]=H]\,,
$$
where $H$ is a chosen to be a graph with degree sequence $\dv_0$. Note that the number of random edges outside $H$ to be exposed is 
${k\choose 2}-{i\choose 2}$, and each of them appears independently and with probability $1/2$. Therefore, the above probability equals to 
$N(\dv_0)/2^{{k\choose 2}-{i\choose 2}}$.
By Corollary~\ref{c:Nd0}, this is asymptotic to  $G(\dv_1)/2^{k\choose 2}=p_k$. 
\qed

\section{Concluding remarks}
Our technique for proving Proposition~\ref{p:ratio} is a rather complicated comparison argument somewhat related to the method of 
switchings used for graphs of similar densities in~\cite{KSVW}. One might be tempted to try proving the result for 
$|S|=\sum_{j=1}^i s_j=2$, as sketched 
in the first part of the proof, and then applying this repeatedly, as in the proof of Corollary~\ref{c:all}, to go from one degree 
sequence to another. However, this seems to provide insufficient accuracy. Similarly, attempts to use switchings directly were not 
successful.

Of independent interest is the following estimate for the probability that a regular graph contains a given subgraph, which gives an 
asymptotic formula provided the sum of the absolute values of $\delta_j=d_j-(i-1)/2$ is a bounded multiple of $k^{3/4}$.
\begin{theorem}\lab{c:probformula}
Assume $i=o(\sqrt \kk)$, with $\kk$ odd. Let $H$ be  a graph on vertex set $[i]$ with degree sequence  $\dv$.
Then the probability that a random $\frac12 (k-1)$-regular graph $G$ on vertex set $[\kk]$ has the induced subgraph $G[i]$ equal to $H$ is
$$
2^{-{i\choose 2}} \exp\left(  \frac{2}{k-i}\sum_{j=1}^i -
\delta_j^2\right)  \exp\left(o(k^{-3/4})\sum_{j=1}^i
|\delta_i|\right)\ .
$$
\end{theorem}
\proof
We may use the argument in the proof of Corollary~\ref{c:all} to jump from $\dv$ to $\dv_0$, using at most $k^{-3/4}\sum_{j=1}^i |\delta_i|$ applications of Corollary~\ref{c:expratio}, and  then apply Corollary~\ref{c:Nd0}. \qed

\end{document}